\documentclass[12pt,a4paper]{amsart}
\usepackage{euscript,amsfonts,amssymb,amsmath,amscd}

\usepackage[T2A]{fontenc}

\usepackage[utf8]{inputenc}        

\usepackage[english,russian]{babel}

\usepackage{hyperref}

\input{epsf}
\newcommand{\op}{\operatorname}
\newcommand{\m}{\mathbb}

\theoremstyle{plain}
\newtheorem{theorem}{Theorem}
\newtheorem{lemma}{Lemma}

\theoremstyle{definition}

\begin{document}

\title {Universality of actions on $\m HP^2$.}

\author{Andrey Kustarev}

\renewcommand{\abstractname}{Abstract}
\renewcommand{\refname}{References}

\begin{abstract}
We show that any eight-dimensional oriented manifold~$M$
possessing a smooth circle action with exactly three fixed points
has the same weights as a certain circle action on~$\m HP^2$. It
follows that the Pontryagin numbers and the equivariant cohomology
of~$M$ coincide with those of~$\m HP^2$. If $M$ admits a cellular
decomposition with three cells, then it is diffeomorphic to~$\m
HP^2$.
\end{abstract}

\maketitle

It is well-known that the quaternionic projective plane $\m HP^2$
admits a circle action with exactly three fixed points. This
action is of particular interest by several reasons. The
quaternionic projective plane is an even-dimensional homogeneous
space admitting no complex or almost complex structure, but
nevertheless having a vast symmetry group. Examples of circle
actions with exactly three isolated fixed points are surprisingly
rare. Known examples of homogeneous spaces $G/H$ admitting an
action of one-dimensional subgroup of $G$ with exactly three fixed
points include $\m CP^2, \m HP^2$ and $\m CaP^2$. Speaking of $\m
HP^2$, one can notice that this manifold also admits a cellular
decomposition with exactly three cells. Considering the standard
circle action on $\m CP^2$ with isolated fixed points, we know
that it is Hamiltonian with respect to the standard symplectic
structure. Since any Hamiltonian is always a Morse function, it
produces  a cellular decomposition of $\m CP^2$ with the cells
corresponding to the fixed points of the original action. In this
respect, $\m HP^2$ behaves similarly to $\m CP^2$, but clearly $\m
HP^2$ has no symplectic structure (since it cannot be almost
complex) and therefore has no Hamiltonian actions. The
connection between circle actions and topology of the manifold in
this case would be more subtle.



We recall that the quaternionic projective plane $\m HP^2$ is
defined as~the~quotient space ($\m H^3_*\,/ \sim$), where $\m H =
\{a + bi + cj + dk\}$, $\m H^3_* = \m H^3 \setminus \{0, 0, 0\}$
and $(x_1, x_2, x_3)\sim (y_1, y_2, y_3)$ iff $(x_1, x_2, x_3) =
(y_1q, y_2q, y_3q)$ for some $q\in \m H\setminus\{0\}$. We view
the circle $S^1$ as~the multiplicative~group of complex numbers
with absolute value~$1$, together with a natural inclusion into
$\m H\setminus\{0\}$ given by the formula $(a + b\cdot i) \to (a +
b\cdot i + 0\cdot j + 0\cdot k)$.

Consider now two families of circle actions on $\m HP^2$:
\begin{enumerate}
{\item The standard action of the form \\ $e^{2\pi it}\cdot
(x_1:x_2:x_3) = (e^{2 k_1\pi it}x_1 : e^{2 k_2\pi it}x_2 : e^{2k_3
\pi it}x_3)$,} {\item The ``semi-integer'' action of the form \\
$e^{2\pi it}\cdot (x_1:x_2:x_3) = (e^{(1+2k_1)\pi it}x_1 :
e^{(1+2k_2)\pi it}x_2 : e^{(1+2k_3) \pi it}x_3)$, \\ $k_1, k_2,
k_3 \in \m Z.$}
\end{enumerate}

In the case of $\m CP^2$ these two families are equivalent, since
the multiplication in $\m C$ is commutative. However, they are
different on $\m HP^2$ -- as we will see later, their sets of
weights at the fixed points are different.

If the circle acts smoothly on an oriented manifold with only
isolated fixed points, then the tangent representation of the
circle is defined at each fixed point. This representation splits
into a sum of two-dimensional irreducible circle representations,
and every real irreducible representation of the circle is
characterized by a single invariant: a positive integer called the
weight of the action. By definition, the sign of a fixed point of
a circle action on an oriented manifold is equal to $(+1)$ if the
orientation of the manifold coincides with the orientation of the
direct sum of the induced two-dimensional representations at this
point and is equal to $(-1)$ otherwise.

\begin{theorem}
\label{T1} The set of weights of any smooth circle action on a
$8$-dimensional closed oriented manifold $M$ with exactly three
fixed points is equivalent to the set of weights of some circle
action on $\m HP^2$ belonging to one of the two families described
above.
\end{theorem}

\begin{theorem}
\label{T2} If a manifold $M$ admits a cellular decomposition with
exactly three cells (or, equivalently, admits a Morse function
with exactly three critical points) in addition to the conditions
of Theorem~1, then $M$ is diffeomorphic to~$\m HP^2$.
\end{theorem}

Note that Theorem \ref{T1} implies that the Pontryagin numbers and
the equivariant cohomology of the manifold coincide with those of
$\m HP^2$. This follows from the localization theorem, whose
precise formulation is given later.

Let us show how one can deduce Theorem \ref{T2} from
Theorem~\ref{T1}. As shown in~\cite{eellskuiper}, any
$8$-dimensional smooth manifold $M$ admitting a cellular
decomposition with exactly three cells is determined up to
diffeomorphism by its Pontryagin numbers. By Theorem~\ref{T1}, the
weights of the circle action 
on $M$ coincide with the weights of some circle action on~$\m
HP^2$. Thus the Pontryagin numbers of $M$ and $\m HP^2$ must also
coincide by the localization theorem, which implies
Theorem~$\ref{T2}$.

The proof of theorem $\ref{T1}$ is the subject of this paper.
Without loss of generality, we may assume the circle action on $M$
to be faithful; otherwise we just factor out the finite subgroup
acting trivially on~$M$. In the case of a faithful action, the GCD
of all weights at each single fixed point is always~$1$.

The main three technical lemmas used in the proof are given
below. In the statements of these lemmas, we consider the general
case of a circle acting smoothly on a $2n$-dimensional oriented
closed manifold $M^{2n}$ with $m$ fixed points. Denote by
$\{a_{ki}\}$ the weights at the fixed point $q_i$, $i=1\ldots m$,
$k=1\ldots n$.

\begin{lemma}[Localization Theorem, \cite{ab}, \cite{toricgenera}]\label{L1}
Denote by $p_{\sigma}(M^{2n})$ the Pontryagin class corresponding
to a symmetric polynomial $\sigma$ in $n$ variables such that
$\deg\sigma\leqslant\frac{n}2$. Then
$$
\bigl(p_{\sigma}(M^{2n}), [M^{2n}]\bigr) = \sum\limits_{i = 1}^m
\op{sign}(q_i)\frac{\sigma(a_{1i}^2,\ldots,a_{ni}^2)}{\prod\limits_{k
= 1}^n a_{ki}}.
$$
\end{lemma}
As usual, we set $(p_{\sigma}(M^{2n}), [M^{2n}]) = 0$ if $\op{deg}
\sigma < n / 2$. In particular, in the case of an action of $S^1$
on an eight-dimensional manifold $M$ with exactly three fixed
points, we obtain two homogeneous equations on the weights of the
action, for the unit class $1\in H^0(M,\m Z)$ and for the first
Pontryagin class $p_1(M)\in H^4(M, \m Z)$:
$$
\frac{\op{sign}(q_1)}{a_{11}a_{21}a_{31}a_{41}} +
\frac{\op{sign}(q_2)}{a_{12}a_{22}a_{32}a_{42}} +
\frac{\op{sign}(q_3)}{a_{13}a_{23}a_{33}a_{43}} = 0,
$$
\begin{multline*}
\op{sign}(q_1)\frac{(a_{11}^2 + a_{21}^2 + a_{31} ^2 + a_{41}^2)}{a_{11}a_{21}a_{31}a_{41}} +\\
+\op{sign}(q_2)\frac{(a_{12}^2 + a_{22}^2 + a_{32} ^2 + a_{42}^2)}{a_{12}a_{22}a_{32}a_{42}} +\\
+\op{sign}(q_3)\frac{(a_{13}^2 + a_{23}^2 + a_{33} ^2 + a_{43}^2)}{a_{13}a_{23}a_{33}a_{43}} = 0.
\end{multline*}

Note that the localization formula for the unit class immediately
implies non-existence of circle actions on closed manifolds having
only one fixed point.

\begin{lemma}
[weights of finite subgroup actions, see \cite{bredon}]\label{L2}
Let $a\ne 1$ be one of the weights at the fixed point $q_i$.
\begin{enumerate}

{\item The set $M^{\m Z_a}$ of points which are invariant under
the action of the subgroup $\m Z_a\subset S^1$ on~$M^{2n}$ is a
disjoint union of smooth connected submanifolds (of possibly
different dimensions) in $M^{2n}$.}

{\item If $a$ is also present among the weights at some fixed
point~$q_j$ lying in the same connected component of $M^{\m Z_a}$
as $q_i$, then the sets of weights at the points $q_i$ and $q_j$
coincide after factorization $\m Z_{>0}\to  \m Z_a/\pm1$.}
\end{enumerate}
\end{lemma}

Lemma \ref{L2} has an interesting implication in the case of a
circle action with exactly three fixed points. Every component of
$M^{\m Z_a}\subset M^{2n}$ is also a manifold equipped with circle
action having only isolated fixed points and, as mentioned above,
there are no components with exactly one fixed point. So if
$mult_j(a)$ is the number of weights at a point $q_j$ which are
multiples of some integer number $a > 1$, then all of the nonzero
numbers $mult_j(a)$ for $j = 1 \ldots 3$ are equal to each other.

\begin{lemma}[weight pairing and weight graph, \cite{musin},
\cite{krichever}]\label{L3} The weights of a~circle action with
isolated fixed points on an oriented manifold can be split into
pairs of equals corresponding to different fixed points.
\end{lemma}

Lemma \ref{L3} implies that the weights at the fixed points $q_1,
q_2, q_3$ on our eight-dimensional manifold $M$ can  be split in
pairs in the following way:
\begin{itemize}
{\item the weights at $q_1$ have the form $a_1, a_2, b_1, b_2$;}
{\item the weights at $q_2$ have the form $a_1, a_2, c_1, c_2$;}
{\item the weights at $q_3$ have the form $b_1, b_2, c_1, c_2$.}
\end{itemize}

Recall that $a_1, a_2, b_1, b_2,c_1, c_2$ are positive integers.
We may assume without loss of generality that $a_1$ is the maximal weight (there
may be other weights equal to it). Applying the localization
theorem for the unit Pontryagin class, we obtain $a_1 > 1$.

Applying Lemma \ref{L2} in our situation, we obtain that there are
exactly three possibilities for the weights $c_1$ and $c_2$:
\begin{enumerate}
{\item $\{c_1, c_2\} = \{a_1 - b_1, a_1 - b_2\};$}
{\item $\{c_1, c_2\} = \{b_1, b_2\}$;}
{\item $\{c_1, c_2\} = \{a_1 - b_1, b_2\}$ and $b_2\ne a_1/2$.}
\end{enumerate}

In any of the three cases above, there is a pair of fixed points
with the opposite signs; otherwise we get a contradiction with the
localization theorem for the unit Pontryagin class.

We will show below that only the first case is possible. For given
numbers $a_1, a_2, b_1, b_2$ and $c_i$ = $a_1 - b_i$, $i=1, 2$,
one can explicitly write out an action on $\m HP^2$ with the
required weight pattern.

\section{Case One}

In this case we have $c_1 = a_1 - b_1, c_2 = a_1 - b_2$. We first
show that the sign distribution of the form $(+ + -)$ for the
points $q_1,q_2,q_3$ is impossible. The remaining two
possibilities are in fact symmetric; we will show that they both
correspond to circle actions on $\m HP^2$ belonging to one of the
two families considered at the beginning of this paper. Namely, we
will show that the action is standard (resp. semi-integer) if and
only if the sum $(a_1+a_2)$ or, equivalently, the sum $b_1+b_2$,
is even (resp. odd).

Suppose that the signs of the points $(q_1, q_2, q_3)$ are $(+ +
-)$. Then we can come to a contradiction by applying the
localization theorem for the Pontryagin classes $1$ and $p_1(M)$.

We obtain the following equations on the weights:
$$
(a_1-b_1)(a_1-b_2) + b_1 b_2 - a_1 a_2 =0,
$$
\begin{multline*}
(a_1-b_1)(a_1-b_2)(a_1^2+a_2^2+b_1^2+b_2^2) +\\
+b_1b_2 ((a_1-b_1)^2 + (a_1-b_2)^2 + a_1^2+a_2^2)-\\
- a_1 a_2 ((a_1-b_1)^2 + (a_1-b_2)^2 + b_1^2 + b_2^2) = 0.
\end{multline*}

The first equation is equivalent to $a_1(a_1-a_2-b_1-b_2) =
-2b_1b_2$. Now we subtract the first equation multiplied by
$(a_1^2+a_2^2+b_1^2+b_2^2)$ from the second one. We obtain
$$
b_1b_2(2a_1^2-2a_1b_1 - 2a_1b_2) - a_1a_2 (2a_1^2 - 2a_1 b_1 -2 a_1 b_2 + b_1^2 + b_2^2 - a_1^2 - a_2^2) = 0.
$$
From the first equation we get $a_1^2-a_1b_1 - a_1b_2 = a_1 a_2 -
2b_1 b_2$. Substituting this into the equation above we get
$$
b_1b_2(2a_1 a_2 - 4b_1 b_2) = a_1 a_2 (2a_1 a_2 - 4b_1 b_2  + b_1^2 +b_2^2 -a_1^2 - a_2^2),
$$
$$
b_1b_2(2a_1 a_2 - 4b_1b_2) = a_1a_2 ( -6b_1b_2 + (b_1+b_2)^2 - (a_1-a_2)^2),
$$
$$
b_1b_2(2a_1 a_2 - 4b_1b_2) = a_1a_2 ( -6b_1b_2 + (b_1+b_2-a_1+a_2)(b_1+b_2+a_1-a_2))
$$
Now recall that $b_1+b_2-a_1+a_2 = 2b_1b_2/a_1$, so we may divide both parts by $2b_1b_2$:
$$
a_1a_2-2b_1b_2 = a_1 a_2 ( -3+(b_1+b_2+a_1-a_2)/a_1),
$$
$$
-2b_1b_2 = a_2 (-3a_1 + b_1 + b_2 - a_2)
$$
Applying again $-2b_1b_2 = a_1(a_1 - a_2 - b_1 - b_2)$, we finally get
$$
a_1^2-a_1a_2-a_1b_1-a_1b_2 = -3a_1a_2 +a_2b_1 + a_2b_2 - a_2^2,
$$
and so $(a_1+a_2)^2 = (a_1 + a_2) \cdot (b_1+b_2)$ and $a_1 + a_2
= b_1 + b_2$. Substituting this into the first equation
$a_1(a_1-a_2-b_1-b_2) = -2b_1b_2$, we also get $b_1b_2=a_1a_2$.
But this implies that one of the weights $b_1,b_2$ is equal to
$a_1$ and so one of the weights $(a_1-b_1)$ or $(a_1-b_2)$ is
zero. This is impossible since the action has only isolated fixed
points.

Now suppose that the signs are of the form $(- + +)$. Then the
localization equation for the unit Pontryagin class is
$$
-(a_1-b_1)(a_1-b_2) + b_1 b_2 + a_1 a_2 =0,
$$
$$
a_1(a_1-b_1-b_2-a_2)=0,
$$
so we have $a_1 = b_1 + b_2 + a_2$.

Consider the circle action on $\m HP^2$ given by the formula
$$
e^{2\pi i t}\cdot (x_1 : x_2 : x_3) = (e^{2\pi it p_1}x_1: e^{2\pi it p_2}x_2 : e^{2\pi it p_3}x_3),
$$
where all of the $p_1,p_2,p_3$ are either integers or
semi-integers.

Direct calculations show that
\begin{itemize}
{\item the weights at the point $(1:0:0)$ are $|p_2\pm p_1|, |p_3
\pm p_1|$;}

{\item the weights at the point $(0:1:0)$ are $|p_2 \pm p_1|$,
$|p_3\pm p_2|$;}

{\item the weights at the point $(0:0:1)$ are $|p_3 \pm p_2|$,
$|p_3\pm p_1|$.}
\end{itemize}

Returning to the original problem, our task is to find numbers
$p_1, p_2, p_3$ corresponding to the initial values of the weights
$a_1,a_2,b_1,b_2$. If we assume $0 \leqslant p_1 < p_2 < p_3$ and
$b_2\leqslant b_1$, then the answer is
$$
a_1 = p_3+p_2,
$$
$$
a_2 = p_3-p_2,
$$
$$
b_1 = p_2+p_1,
$$
$$
b_2 = p_2-p_1.
$$
Then $a_1 = a_2 + b_1 + b_2$ and the fixed point $q_1$ corresponds
to $(0:1:0)$. The next point containing the weight $a_1$ is
$(0:0:1)$ and it corresponds to $q_3$, so $(1:0:0)$ corresponds to
$q_2$.

Now suppose that the signs are of the form $(+ - +)$. Then the
localization equation is
$$
(a_1-b_1)(a_1-b_2) - b_1 b_2 + a_1 a_2 =0,
$$
$$
a_1(a_1-b_1-b_2+a_2)=0,
$$
so $a_1 + a_2 =  b_1 + b_2$. Parameters of the corresponding
circle action on $\m HP^2$ now may be found as in the previous
case.

\section{Case Two}
Here we suppose that $c_1=b_1, c_2=b_2$.
Applying the localization formula for the unit Pontryagin class we
obtain
$$
\pm\frac1{a_1 a_2 b_1 b_2} \pm \frac1{a_1 a_2 b_1 b_2} \pm \frac1 {b_1 b_2 b_1 b_2} = 0,
$$
and so the signs of the fixed points should have the form $( + +
-)$. So $2 b_1 b_2 = a_1 a_2$. Now we write the localization
formula for the Pontryagin class $p_1(M)$:
$$
2\frac{a_1^2+a_2^2+b_1^2+b_2^2}{a_1 a_2 b_1 b_2} - \frac {b_1^2 + b_2^2 + b_1^2 + b_2^2}{b_1 b_2 b_1 b_2} = 0.
$$
Simplifying this expression, we get
$$
a_1^2 + a_2^2 = b_1^2 + b_2^2.
$$
But $a_1 a_2 = 2 b_1 b_2 > 0$. Subtracting this from the last
equation, we get $(a_1-a_2)^2 < (b_1 - b_2)^2$. Considering
$b_2\leqslant b_1$ and using $a_2\leqslant a_1$ we get that
$a_1-a_2<b_1-b_2$. Similarly, $(a_1+a_2)^2 > (b_1+b_2)^2$, and so
$a_1+a_2
> b_1+b_2$. Summing up these inequalities implies $a_2>b_2$. Since
$a_1\geqslant b_1$, the sums of squares $(a_1^2 + a_2^2)$ and
$(b_1^2 + b_2^2)$ are not equal. This is a contradiction.

\section{Case Three}

Finally, suppose that $c_1 = a_1 - b_1, c_2 = b_2$. We will see
that this is possible only when $b_2 = a_1/2$ (so $c_i = a_1 -
b_i$, just as in the Case~$1$). To show that we will have to check
out all three possible cases of the fixed point sign
distributions. So the proof splits here in three more cases.

Now the localization equations have the following form:
$$
\pm (a_1-b_1)b_2 \pm b_1 b_2 \pm a_1 a_2 = 0,
$$
\newpage
\begin{multline*}
\pm (a_1-b_1)b_2(a_1^2+a_2^2+b_1^2+b_2^2)\pm\\
\pm b_1b_2 ((a_1-b_1)^2+b_2^2+a_1^2+a_2^2)\pm\\
\pm a_1a_2 ((a_1-b_1)^2+b_2^2+b_1^2+b_2^2)=0.
\end{multline*}

Here we denote by $\pm$ the sign of the corresponding fixed point
$q_1,q_2,q_3$; recall that there are exactly two pluses and one
minus among these signs.

Suppose that the signs of points $q_1,q_2,q_3$ are $(+ + -)$,
respectively. Then the localization equation for the unit
Pontryagin class is
$$
(a_1-b_1)b_2 + b_1b_2 - a_1a_2=0,
$$
so $a_2 = b_2 = d$ (say).

Now consider the localization equation for the Pontryagin class
$p_1(M)$ and subtract from it the equation for the unit Pontryagin
class multiplied by $(a_1^2+d^2+b_1^2+d^2)$. We get
$$
b_1d(a_1^2-2a_1b_1) - a_1d(-2a_1b_1+b_1^2)=0.
$$
After simplification we obtain that $a_1=b_1$, so the weight $c_1
= a_1 - b_1$ is zero. This is impossible since the action has only
isolated fixed points.

Now consider the case of sign distribution $(+ - +)$. Then
$$
b_2(a_1-b_1)-b_1b_2+a_1a_2=0
$$
and $b_2(2b_1-a_1)=a_1 a_2$. Again, considering the localization
equation for the Pontryagin class $p_1(M)$ and subtracting from it
the equation for the unit Pontryagin class multiplied by
$(a_1^2+a_2^2+b_1^2+b_2^2)$, we get
$$
-b_1b_2(a_1^2-2a_1b_1) + a_1a_2(b_1^2-2a_1b_1+b_2^2-a_2^2) = 0.
$$
But $-b_2(a_1-2b_1) = a_1a_2$, so we can make the corresponding substitution in the equation and drop $a_1a_2$. The result is
$$
b_1^2 - a_1b_1 + b_2^2 -a_2^2=0,\quad (b_2-a_2)(b_2+a_2)=(a_1-b_1)b_1.
$$

We know that $a_1-b_1>0$ and so $b_2-a_2>0$ and $b_2 > 1$. The
weight $b_2$ is present at all three fixed points. Consider now
the connected component of $M^{\m Z_{b_2}}$ containing the points
$q_1,q_2,q_3$. Each of these points contains at least one weight
that is a multiple of $b_2$, so the numbers of such weights at
each of $q_1, q_2, q_3$ should be equal by Lemma~\ref{L2}. Note
that there are at least two weights equal to $b_2$ at the point
$q_3$.

Since $a_2<b_2$ and the numbers of weights that are multiples of
$b_2$ coincides at $q_1$ and $q_2$, the weights $b_1$ and
$a_1-b_1$ are either both multiples of~$b_2$, or none of them is a
multiple of~$b_2$. If they are both multiples of $b_2$, then all
weights at $q_3$ are also multiples of $b_2$ and since
$b_2>a_2\geqslant 1$, the circle action is not faithful.

So the only remaining possibility is that $a_1=kb_2$. One can now
write down the localization formula for the unit Pontryagin class
and the action of the quotient circle $S^1/{\m Z_{b_2}}$ on the
four-dimensional component $M^{\m Z_{b_2}}$ containing
$q_1,q_2,q_3$. A simple explicit calculation shows that $k$=$2$ is
the only possible value. So $b_2=a_1/2$ in this case.

Finally, we consider the sign distribution $(- + +)$. Then the
localization equation for the unit Pontryagin class has the form
$$
-(a_1-b_1)b_2+b_1b_2+a_1a_2=0, \quad a_1a_2  = (a_1 - 2b_1)b_2,
$$
and the localization equation for the Pontryagin class $p_1(M)$ yields
$$
b_1b_2(a_1^2-2a_1b_1) + a_1a_2(b_1^2-2a_1b_1+b_2^2-a_2^2) = 0,
$$
$$
a_1 a_2 b_1 a_1+ a_1a_2(b_1^2 - 2a_1b_1 + b_2^2 - a_2^2)=0,
$$
$$
b_1^2-a_1b_1+b_2^2-a_2^2=0,\quad (b_2-a_2)(b_2+a_2)=(a_1-b_1)b_1.
$$
Since $b_2>a_2$, we may repeat all steps as above and get that $b_2=a_1/2$. This finishes the proof.


\end{document}